\documentclass{article}

\usepackage{arxiv}

\usepackage[utf8]{inputenc} 
\usepackage[T1]{fontenc}    

\usepackage{booktabs}       
\usepackage{nicefrac}       
\usepackage{microtype}      
\usepackage{lipsum}		
\usepackage{doi}

\usepackage{hyperref, url}
\usepackage{amsfonts, amsmath, amssymb, parskip, bbm}
\usepackage{xcolor, graphicx, tikz-cd, caption, subcaption}
\usepackage{algorithmic, algorithm}
\usepackage{verbatim}
\usepackage{mathrsfs}

\usepackage{amsthm}  
\newtheorem{theorem}{Theorem}
\newtheorem{remark}{Remark}
\newtheorem{lemma}{Lemma}

\newcommand{\NN}{{\mathtt N}}
\newcommand{\DD}{{\mathtt D}}
\newcommand{\RR}{{\mathtt R}}
\newcommand{\PP}{{\mathtt P}}

\title{Rational function approximation with normalized positive denominators}

\author{
  \hspace{1mm}James Chok \\
  School of Mathematics and \\
  Maxwell Institute for Mathematical Sciences, \\
  The University of Edinburgh,\\
  Edinburgh,\\
  EH9 3FD, \\
  United Kingdom\\
  \texttt{james.chok@ed.ac.uk} \\
    \And
  Geoffrey M. Vasil \\
  School of Mathematics and \\
  Maxwell Institute for Mathematical Sciences, \\
  The University of Edinburgh,\\
  Edinburgh,\\
  EH9 3FD, \\
  United Kingdom\\
  \texttt{gvasil@ed.ac.uk} \\
}

\hypersetup{
pdftitle={Rational function approximation with normalized positive denominators},
pdfsubject={math.NA},
pdfauthor={James Chok, Geoffrey M. Vasil},
pdfkeywords={Rational approximation, Suprious poles, Spectral methods, Data analysis, Derivative penalty},
}

\begin{document}
\maketitle

\begin{abstract}
Recent years have witnessed the introduction and development of extremely fast rational function algorithms. Many ideas in this realm arose from polynomial-based linear-algebraic algorithms. However, polynomial approximation is occasionally ill-suited to specific challenging tasks arising in several situations. Some occasions require maximal efficiency in the number of encoding parameters whilst retaining the renowned accuracy of polynomial-based approximation. One application comes from promoting empirical pointwise functions to sparse matrix operators. Rational function approximations provide a simple but flexible alternative (actually a superset), allowing one to capture complex non-linearities. However, these come with extra challenges: i) coping with singularities and near singularities arising from a vanishing denominator, and ii) a non-uniqueness owing to a simultaneous renormalization of both numerator and denominator. 

We, therefore, introduce a new rational function framework using manifestly positive and normalized Bernstein polynomials for the denominator and any traditional polynomial basis (e.g., Chebyshev) for the numerator. While an expressly non-singular approximation slightly reduces the maximum degree of compression, it keeps all the benefits of rational functions while maintaining the flexibility and robustness of polynomials. We illustrate the relevant aspects of this approach with a series of derivations and computational examples. 
\end{abstract}


Accepted to SIAM Journal on Scientific Computing

\section{Introduction}\label{sec:intro}
Function approximation involves approximating a target function $f:\mathbb{R}\to\mathbb{R}$ with a `simpler' function $g:\mathbb{R}\to\mathbb{R}$. This can be done for reasons including machine learning \cite{ml_1, ml_3, ml_4, ml_2} and numerical computations of transcendental functions \cite{transcendental_1, transcendental_2}. This paper considers rational polynomial approximations
\begin{equation}
    f(x)\approx r(x) = \frac{p(x)}{q(x)}, \quad\text{where}\quad p(x)\in\mathcal{P}(n)\ \text{and}\ q(x)\in\mathcal{P}(m),
\end{equation}
where $\mathcal{P}(n)=\{\text{Degree } n \text{ polynomials}\}$ and we denote $\mathcal{R}(n, m) = \{p(x)/q(x)\ |\ p(x)\in\mathcal{P}(n)\ \text{and}\ q(x)\in\mathcal{P}(m)\}$ as the space of $(n,m)$ rational polynomials. 

{A $\mathcal{R}(n, n)$ rational polynomial often converges at the same rate as a $\mathcal{P}(2n)$ degree polynomial \cite{Trefethen2023}. This factor of 2 improvement is important in our motivating application of solving non-constant-coefficient differential equations}
\begin{equation}
\sum_{k=0}^{K} c_{k}(x) \, D^{k} u(x) \ = \ f(x).
\end{equation}
where we assume $f(x) \in \mathcal{P}(n)$ for some $n$, {and $c_k(x)$ are smooth functions. We do not consider the case in which the coefficients are non-smooth.}

Chebyshev spectral methods are a common approach to equations of this form \cite{chebyshev_spectral_methods, chebop, chebyshev_collocation, chebysev_spectral_radiative_transfer, spectral_methods_matlab}. In that case, the solution is expanded in terms of polynomial basis functions, and the derivatives map sparsely to an alternative basis set.
The catch is that we must expand the non-constant coefficients, $c_{k}(x)$, in a polynomial basis and promote the result to matrix operations. Such an approach works especially well in one-dimension \cite{spectral_methods_sparse, dedalus, vasil_disk, ball1, ball2, papadopoulos2023building} where the typical degree of a coefficient, $c_{k}(x)$, is much less than the requirements of the solution, e.g., $\sim 10$ versus $\sim 1000$ respectively.  {The downside is that the final overall matrix representation becomes dense with each additional polynomial degree ever more required by the coefficients; the matrix bandwidth scales with the degree.} The cubic scaling of operator bandwidth throttles the overall algorithm. The problem is even more stark in multiple dimensions \cite{OTV2019, OTV2020, Ellison_Julien_Vasil_2022, Olver2020Fast, burns2022corner, Townsend2015Automatic}.

{
Instead, we could represent the coefficients as 
\begin{equation}
c_{k}(x) \ = \ \frac{p_{k,N}(x)}{q_{M}(x)}, \quad \text{where} \quad p_{k,N}(x) \in \mathcal{P}(N), \quad q_{M}(x) \in \mathcal{P}(M).
\end{equation}
Suppose all the coefficients have a common denominator $q_{M}(x)$, independent of $k$. 
We can thus use a smaller maximum degree $N$ than a pure polynomial representation requires. Clearing the denominator allows us to solve more sparsely 
\begin{equation}
\sum_{k=0}^{K} p_{k,N}(x) \, D^{k} u(x) \ = \ q_{M}(x) \, f(x).
\end{equation} 
Particularly, the factor of 2 improvement suggests that the operator bandwidth is reduced by half. We discuss applications like this and more in our examples section. }

{
However, a prevalent problem in rational approximations is the production of spurious poles  \cite{berrut_poleless_barycentric, barycentric_poleless, rational_fragile}. This occurs when $q_{M}(x)$ has roots in or near the relevant part of the domain.}

To that end, we propose a new \textit{Bernstein Denominator Algorithm}, where we force $q(x) > 0$ in the approximation interval by using Bernstein polynomials with normalized and positive coefficients. This provides a robust method to produce rational approximations without poles in the approximation domain. This method works well on noisy data, consistently producing approximations without poles while other rational approximations do. 

\section{Previous Works}
A simple way of representing rational polynomials in $\mathcal{R}(n, m)$ is in the form
\begin{eqnarray}\label{eq:normalized_rational_polynomial}
    r(x) = \frac{p(x)}{q(x)} = \frac{a_0 + a_1 x + \cdots + a_n x^n}{1 + b_1 x + \cdots + b_m x^m},
\end{eqnarray}
where the one in the denominator represents a normalizing factor; otherwise, it is invariant under rescalings of the numerator and denominator, leading to non-unique solutions.

Finding the coefficients can be framed as a least-squares problem, 
\begin{eqnarray}\label{eq:og_rational_problem}
    \min_{\substack{a_0,\ldots, a_n\\ b_1,\ldots, b_m}}\sum_i \epsilon_i^2, \quad \text{where} \quad \epsilon_i = \frac{p(x_i)}{q(x_i)} - f(x_i),
\end{eqnarray}
where $\{x_i\}$ are distinct points over $\mathbbm{R}$. This, however, proves hard to solve as it becomes a nonlinear optimization problem. Instead, one typically considers the linearized residuals to which the solution can be found via singular value decomposition (SVD) \cite{approximation_theory_and_practice},
\begin{eqnarray}
    \min_{\substack{a_0,\ldots, a_n\\ b_1,\ldots, b_m}}\sum_i \varepsilon_i^2, \quad \text{where} \quad \varepsilon_i = p(x_i) - f(x_i)q(x_i).
\end{eqnarray}

\subsection{SK Algorithm}
{The SK Algorithm} was introduced by Sanathanan and Koerner \cite{sk_algorithm} to reduce the deficiencies present in the formulation of the linearized residuals. They note that the linearized residuals $\varepsilon_i$ would result in bad fits for lower frequency values and larger errors when $q(x)$ has poles in the complex plane. Instead, they consider an iterative procedure in which, at iteration $t$, one minimizes the reweighted residuals
\begin{eqnarray}
    \min_{\substack{a_0,\ldots, a_n\\ b_1,\ldots, b_m}}\sum_i \mathcal{E}_{t,i} ^2, \quad \text{where} \quad \mathcal{E}_{t, i} = \frac{\varepsilon_{t, i}}{q_{t-1}(x_i)} = \frac{p_t(x_i) - f(x_i) q_t(x_i)}{q_{t-1}(x_i)},
\end{eqnarray}
which provides a better reflection of the original rational approximation problem (\ref{eq:og_rational_problem}) \cite{stablized_sk_algorithm}. Taking $q_0(x)=1$ for the first iteration, each iteration can similarly be solved using SVD \cite{stablized_sk_algorithm}.

\subsection{AAA Algorithm}
{The AAA algorithm}, by Nakatsukasa, S\`ete, and Trefethen \cite{AAA}, is a foundational work that set off a vast amount of research \cite{barycentric_periodic_functions, barycentric_minmax, aaa_eigenvalue_problems, AAA_Lawson, p_aaa} and is based on rational Barycentric forms
\begin{eqnarray}\label{eq:barycentric_form}
    r(x) = \frac{n(x)}{d(x)} =  \sum_{k=1}^n\frac{w_k f_k}{x - x_k} \bigg{/}\sum_{k=1}^n \frac{w_k}{x - x_k},
\end{eqnarray}
with $f_k\in\mathbb{C}$ and $w_k \in \mathbb{C} - \{0\}$. As shown in \cite{AAA}, these forms range over the set of $\mathcal{R}(n - 1, n - 1)$ that have no poles at the points $x_k$. Moreover, they have the property that $r(x_k) = f_k$ for all $k$. 

First setting $f_k = f(x_k)$, the AAA algorithm works by finding weights that minimize the residuals 
\begin{eqnarray}
    &&\min_{w\in S^{J-1}} \sum_{i}\varepsilon_i^2,\quad \text{where}\quad S^{J-1} = \{w\in\mathbb{R}^{J}\ |\ \|w\| = 1\},\\
    &&\text{and}\quad \varepsilon_i = n(z_i) - f(z_i)d(z_i)= \sum_{k=1}^{J}\frac{w_k f_k}{z_i - x_k} - f(z_i)\sum_{k=1}^{J}\frac{w_k}{z_i - x_k},
\end{eqnarray}
on a set of distinct support points $\{z_i\}_i$ with $z_i \neq x_k$ for all $i, k$. When $|\{z_i\}_i| \geq J$, a unique minimizer can be found using SVD \cite{AAA}.

\subsection{Barycentric forms with no poles}
{Berrut \cite{berrut_poleless_barycentric} first studied barycentric rational polynomials with no poles in the domain and later generalized by Floater and Hormann \cite{barycentric_poleless}. They take the form
\begin{eqnarray}
    r(x)\ =\ \frac{\sum_{i=1}^{n-d}\lambda_i(x)p_i(x)}{\sum_{i=1}^{n-d}\lambda_i(x)},\quad\text{where}\quad \lambda_i(x)\ =\ \frac{(-1)^i}{(x-x_i)\cdots (x-x_{i+d})},
\end{eqnarray}
where $p_i\in\mathcal{P}(d+1)$ that interpolates $f$ at $x_i,\ldots, x_{i+d}$, and $0\leq d\leq n$. Berrut rational polynomial is the special case when $d=0$, reducing $p_i(x)=x_i$.}

{This can be simplified to the familiar Barycentric form \cite{barycentric_poleless}
\begin{eqnarray}
    r(x) = \sum_{k=0}^n\frac{w_k\, f_k}{x-x_k}\bigg{/}\sum_{k=0}^{n}\frac{w_k}{x-x_k},\ \text{where}\ \   w_k = \sum_{i\in J_k}(-1)^i\sum_{j=1,j\neq k}^{i+d}\frac{1}{x_k-x_j},
\end{eqnarray}
and $J_k=\{i\in[0,\ldots, n] | k - d\leq i\leq k\}$. In this form, it is easy to see that $r\in \mathcal{R}(n,n)$. }
{An important distinction to the AAA is that AAA uses the linearized residuals to find the optimal values of $w_k$. In contrast, Floater and Hormann rational polynomials fix $w_k$ allowing pole-free rational polynomials.}

\section{Rational Bernstein Denominator Algorithm} 
{In this paper, we are interested in representing smooth functions, $f(x)$, that are guaranteed to be pole-free in $[0,1]$ using rational polynomials. We are not interested in representing a larger class of functions. To ensure this, we require the denominator to be strictly positive in $[0,1]$}
\begin{eqnarray}
f(x) \ \approx \ \RR_{N,M}(x) \ = \ \frac{\NN_{N}(x)}{\DD_{M}(x)} \quad \text{with} \quad \DD_{M}(x)\ >\ 0, \quad \forall x\in[0, 1].
\end{eqnarray}

To enforce strict positivity in the denominator, we introduce the family of Bernstein polynomials, named after Sergei Natanovich Bernstein for his probabilistic proof of Weierstrass approximation theorem \cite{weierstrass_theorem_proof}. They are defined as 
\begin{eqnarray}
    \mathcal{B}^{(n)}_{k}(x)\ =\ \binom{n}{k}\, x^{k}\, (1-x)^{n - k},\quad \text{where}\quad 0\leq k\leq n,
\end{eqnarray}
for $x\in[0, 1]$. Noting that $\mathcal{B}_k^{(n)}(x)\geq 0$ in $[0, 1]$, the denominator can  thus be represented as 
\begin{eqnarray}
\DD_{M}(x) \ = \ \sum_{m=0}^{M} w_{m} \, \mathcal{B}_{m}^{(M)}(x),
\end{eqnarray}
where positivity and normalization can be enforced by requiring $w_m\geq 0$ for all $m$, with $\sum_{m=0}^Mw_m=1$. This implies that the weights lie in a probability simplex 
\begin{eqnarray}
\Delta^{M + 1} \ = \ \left\{ w\in\mathbb{R}^{M + 1}\ \bigg{|}\  0\leq w_{m} \leq 1 \ \text{and}\ \sum_{m=0}^{M} w_{m} = 1 \right\}.
\end{eqnarray}

{Note that for all $k$, $\mathcal{B}^{(n)}_k(x)>0$ for $x\in(0,1)$, and $\mathcal{B}^{(n)}_0(0),\mathcal{B}^{(n)}_n(1)>0$. Thus, strict positivity can be enforced by adding the additional constraint $w_0,w_M>0$.}

The numerator is represented as
\begin{eqnarray}
\NN_{N}(x) \ = \ \sum_{n=0}^{N}a_{n}\,  p_{n}(x), \quad \text{where} \quad a_{n} \in \mathbb{R},
\end{eqnarray}
given a suitable family of basis functions (usually orthogonal polynomials). 
No additional constraints on the spectral coefficients $a_n$ are needed, as the denominator is normalized.

Thus, we consider rational functions of the form
\begin{eqnarray}\label{eq:bernstein_rational_function}
    f(x)\ \approx\ \RR_{N,M}(x)\ =\ \sum_{n=0}^N\, a_n\, p_n(x)\ \bigg{/}\ \sum_{m=0}^M\, w_m\, \mathcal{B}_m(x)
\end{eqnarray}
with $w\in\Delta^{M+1}$ and $a\in\mathbb{R}^{N+1}$. We present a method to find these coefficients in Section \ref{sec:residuals}.

\subsection{Lagrange's Polynomial Root Bound}
Another way to guarantee rational approximations without roots in $[0, 1]$ is to use bounds on the roots of polynomials. In particular, for a $n$-degree polynomial, $p(x)=a_0 + a_1x + \cdots + a_n x^n$, Lagrange's Bound states that the magnitudes of its $n$-roots, $\{x_i\}_{1\leq i\leq n}$, can be bounded by
\begin{eqnarray}
    \frac{1}{\max\bigg{\{} 1,\ \sum_{i=1}^n\, \left|{a_i}/{a_0}\right|\bigg{\}}}\ \leq\ |x_i| \quad \text{for}\quad 1\leq i\leq n.
\end{eqnarray}
Since we care about polynomials in the normalized denominator, \textit{i.e.} $a_0=1$, to be free of poles in $[0, 1]$, one can consider polynomials of the form
\begin{eqnarray}\label{eq:lagrange_polynomial}
\widetilde{\mathcal{P}}(n)\ =\ \bigg{\{}1\, +\, a_1x\, +\, \cdots\, +\, a_n x^n, \quad\text{with}\quad \sum_{i=1}^n\,|a_i|\ \leq \ 1\bigg{\}},
\end{eqnarray}
to guarantee no roots inside $[0, 1]$. However, the Bernstein formulation is a superset of this problem. 
{
\begin{lemma}
    For any $p(x)\in \widetilde{\mathcal{P}}(n)$, there exists a $c\in[0,1]$ such that $p(x) = (1 - c) + c\, q(x)$, where $q(x)\in\overline{\mathcal{P}}(n)$, the closure of $\widetilde{\mathcal{P}}(n)$,
    \begin{eqnarray}
    \overline{\mathcal{P}}(n)\ =\ \bigg{\{}1\, +\, a_1x\, +\, \cdots\, +\, a_n x^n, \quad\text{with}\quad \sum_{i=1}^n|a_i|\ =\ 1\bigg{\}}.
\end{eqnarray}
\end{lemma}
}
\begin{proof}
Since $p(x)\in\widetilde{\mathcal{P}}(n)$, it can be written as
\begin{align}
    p(x)\ :=\ 1\, +\, a_1x\, +\, \cdots\, +\, a_n x^n\ &=\ (1-c)\ +\ c\left(1 + \frac{a_1}{c}x +\cdots + \frac{a_n}{c}x^n\right),
\end{align}
where $c=\sum_i|a_i|\leq 1$. It follows that $\sum_i|a_i/c| = 1$, and
\begin{align}
    q(x)\ :=\ 1 + \frac{a_1}{c}x +\cdots + \frac{a_n}{c}x^n\ \in\ \overline{\mathcal{P}}(n).
\end{align}
\end{proof}
{
By properties of the Bernstein polynomials,
\begin{align}
    \sum_{j=0}^n\, \mathcal{B}^{(n)}_j(x)\ =\ 1\quad\text{for all}\quad x\in[0, 1].
\end{align}
Thus, if all polynomials in $\overline{\mathcal{P}}(n)$ can be written in terms of Bernstein basis polynomials with positive coefficients, so can polynomials in $\widetilde{\mathcal{P}}(n)$.}

\begin{theorem}
Any polynomial in $\overline{\mathcal{P}}(n)$ has an equivalent Bernstein basis representation with coefficients that are all positive, i.e., there exist Bernstein coefficients $b_i\geq 0$ such that
\begin{align}
    p(x)\ :=\ 1 + a_1x + \cdots + a_n x^n\ =\ \sum_{i=0}^n\, b_i\, \mathcal{B}^{(n)}_i(x), \quad\text{with}\quad \sum_{i=1}^n\, |a_i|\ =\ 1.
\end{align}
\end{theorem}
\begin{proof}
Since $p\in \overline{\mathcal{P}}(n)$, then $\sum_i|a_i|=1$, and hence
\begin{align}
    p(x)\ &=\ \sum_{a_i\geq 0}\,|a_i|\,(1 + x^i)\ +\ \sum_{a_i < 0}\,|a_i|\,(1 - x^i).
\end{align}
Thus, to account for the sign, these polynomials can be written as a convex combination of
\begin{align}
    M_{i, \pm}(x)\ =\ 1 \pm x^i,\quad \text{for}\quad 1\leq i\leq n,
\end{align}
with positive weight given to one of $M_{i,\pm}$ with the correct sign and zero to the other. That is, there exists positive weights, $w_{i, s}\geq 0$, such that 
\begin{align}
    p(x)\ =\ \sum_{s\in\{+, -\}}\,\sum_{i=1}^n\, w_{i, s}\, M_{i, s}(x)\quad \text{and}\quad \sum_{s\in\{+, -\}}\, \sum_{i=1}^n\, w_{i, s}\ =\ 1.
\end{align}

From the Binomial and Vandermonde's identity, the monomial coefficients, $a_i$, can be expressed in Bernstein coefficients via $b_j = \sum_{k=0}^j a_k \binom{j}{k} / \binom{n}{k}$. Thus
\begin{align}
    M_{i, \pm}(x)\ =\ 1 \pm x^i\ =\ \sum_{j}\, b_j\, \mathcal{B}^{(n)}_j(x)\quad\text{where}\quad b_j\ =\ 1 \pm \frac{\binom{j}{i}}{\binom{n}{i}}\mathbbm{1}_{i\leq j},
\end{align}
and $\mathbbm{1}_{i\leq j} = 1$ if $i\leq j$ and $0$ otherwise. It is clear that $b_j\geq 0$ for $M_{i, +}$. For $M_{i, -}$, since $j\leq n$, it follows from the geometric interpretation of combinations that $\binom{j}{i}\leq\binom{n}{i}$.
Thus $b_j\geq 0$ for $M_{i, \pm}$. Hence, each basis function, $M_{i, \pm}$, represents a positively weighted sum of Bernstein polynomials.

Since every polynomial in $\overline{\mathcal{P}}(n)$ can be written as a convex combination of $M_{i, \pm}$, it follows that it also can be written as a convex combination of Bernstein polynomials with positive coefficients. Thus, every polynomial in $\overline{\mathcal{P}}(n)$ can be written as a sum of Bernstein polynomials with positive coefficients.
\end{proof}

\begin{remark}
The Lagrange bound is not sharp, e.g. $ (\mathcal{B}_0^{(n)} + \mathcal{B}_n^{(n)})/2$ has no roots in $[0, 1]$.
\end{remark}
\begin{proof}
The Bernstein coefficients, $b_i$, can be converted to monomial coefficients, $a_j$, via
\begin{eqnarray}
    a_j\ =\ \sum_{k=0}^j\, (-1)^{j-k}\, \binom{n}{j}\, \binom{j}{k}\, b_k. 
\end{eqnarray}
Explicitly, for $1\leq j\leq n-1$, the monomial coefficients are given by
\begin{align}
    a_0\ =\ \frac{1}{2}, \quad  a_j\ =\ (-1)^j\, \frac{1}{2}\, \binom{n}{j}, \quad\text{and}\quad a_n\ =\ \frac{(-1)^{n}\ + 1}{2}. 
\end{align}
Now,
\begin{align}
    \sum_{j=1}^n\, \bigg{|}\frac{a_j}{a_0}\bigg{|}\ \geq\ \bigg{|}\frac{a_1}{a_0}\bigg{|}\ =\ \binom{n}{1}\ =\ n.
\end{align}
Therefore, Lagrange's lower bound can be bounded by
\begin{align}
    \frac{1}{\max\bigg{\{}1, \sum_{j=1}^n |a_j/a_0|\bigg{\}}}\ \leq\ \frac{1}{n}\ <\ 1.
\end{align}

As such, a normalized-positive Bernstein basis polynomial exists in which Lagrange's lower bound yields a bound on the magnitudes of the roots that is less than one.
\end{proof}

\begin{remark}
While Lagrange's (and Cauchy's) bound yield bounds on roots within a radius of the origin, the Bernstein polynomials offer an alternative way to bound roots on the positive real line. In particular, for a given monomial
\begin{eqnarray}
    p(x)\ =\ a_0\, +\, a_1x\, +\, \cdots \, + a_nx^n,
\end{eqnarray}
if the Bernstein coefficients
\begin{eqnarray}\label{eq:monomials_to_bernstein_coef}
    b_j\ =\ \sum_{k=0}^j\ \frac{\binom{j}{k}}{\binom{n}{k}}\ a_k,
\end{eqnarray}
are all positive for $0\leq j\leq n$, then there are no roots in $(0, 1)$. If $b_0=0$, there is a root at $x=0$ and none otherwise. Similarly, if $b_n=0$, there is a root at $x=1$ and none otherwise.
\end{remark}

\begin{remark}
To guarantee rational fits without poles, one can use a method like the AAA or SK algorithm and check for poles afterward. However, root-finding algorithms can be unstable, giving inaccurate results \cite{numeric_polynomial_roots_unstable, poles_who_cares}.
\end{remark}

{
\begin{theorem}[Weierstrass approximation theorem \cite{weierstrass_theorem_proof}]
    Let $f$ be a continuous function on $[0,1]$, and it's degree $N$ Bernstein representation
    \begin{eqnarray}
        B_N(f)(x)\ =\ \sum_{j=0}^N f(j/N)\, \mathcal{B}_j^{(n)}(x).
    \end{eqnarray}
    Then for any $\varepsilon>0$, there exists an $N$ such that for all $n>N$ and all $x\in[0,1]$
    \begin{eqnarray}
        |B_N(f)(x) - f(x)|\ <\ \varepsilon.
    \end{eqnarray}
\end{theorem}
It follows from the Weierstrass approximation theorem that any strictly positive continuous function $f$ can be approximated by Bernstein polynomials with strictly positive coefficients. However, if $f\in\mathcal{P}(m)$ and strictly positive on $[0,1]$, it can be exactly represented by Bernstein polynomials with strictly positive coefficients. }

{
\begin{theorem}\label{theorem:polynomial_in_bernstein}
    Let $p(x)\in\mathcal{P}(m)$ and strictly positive on $[0,1]$ with
    \begin{eqnarray}
        p(x)\ =\ \sum_{j=0}^n b_j\, \mathcal{B}_j^{(n)}(x)
    \end{eqnarray}
    its corresponding Bernstein polynomial representation for $n\geq m$. For any 
    \begin{eqnarray}
        0\ <\ \varepsilon\ <\ \min_{x\in[0,1]}p(x),
    \end{eqnarray}
    there exists a $N\geq m$ such that, for all $n>N$, $b_j>\varepsilon$ for all $j=0,\ldots, n$.
\end{theorem}
}

\begin{proof}
{
    Let $p(x)=a_0+a_1x + \ldots + a_nx^m$. Since Bernstein polynomials form a polynomial basis, and $n\geq m$, there exists $b_j$ such that
    \begin{eqnarray}
        p(x) \ =\ \sum_{j=0}^n\, b_{j}\, \mathcal{B}_j^{(n)}(x), 
    \end{eqnarray}
    where $b_j$ is defined in (\ref{eq:monomials_to_bernstein_coef}). Using the bounds 
    \begin{eqnarray}
        \left(\frac{n}{k}\right)^k\ \leq\ \binom{n}{k}\ \leq\ \frac{n^k}{k!}\ \implies\ \frac{\binom{j}{k}}{\binom{n}{k}}\ \leq\ \frac{j^k}{k^k}\frac{k!}{n^k}\ \leq\ \frac{j^k}{n^k}
    \end{eqnarray}
    yields
    \begin{eqnarray}
    b_{j}\ &=&\ p(j/n)\ -\ \sum_{k=0}^{\min(j,m)}\left[ \left(\!\frac{j}{n}\!\right)^{\! k}\  -\ \frac{\binom{j}{k}}{\binom{N}{k}}\right] a_{k} \\ 
    &\ge&\ \min_{x\in[0,1]} p(x)\ -\ \max_{k} a_{k}\, \sum_{k=0}^{m}\left[ \left(\!\frac{j}{n}\!\right)^{\! k}\  -\  \frac{\binom{j}{k}}{\binom{n}{k}}\right],\\
    &\geq&\ \min_{x\in[0,1]}p(x)\ -\ \frac{2154}{15\,625} \frac{m^2}{n}\, \max_k a_k, 
    \end{eqnarray}
    where the last inequality comes from Lemma \ref{lemma:combination_bound}.
}

{
    Since $p(x)$ is strictly positive, it follows that $\min_{x\in[0,1]}p(x)>0$. Moreover, since $m$ is fixed, $n$ can be chosen arbitrarily large to make all $b_j$ positive. Taking
    \begin{eqnarray}\label{eq:positive_bernstein_bound}
        N\ =\ \max\left(m, \left \lceil\frac{2154\, m^2}{15\,625}\frac{\max_k a_k}{\min_{x\in[0,1]}p(x)}\right\rceil\right),
    \end{eqnarray}
    where $\lceil\cdot\rceil$ is the ceiling operator, yields the proof
}
\end{proof}

{
\begin{lemma}\label{lemma:combination_bound}
    For $n \ge 1$ and $0 \leq j, m \leq n$,
    \begin{align}
        \sum_{k=0}^{m} \left[ \frac{j^k}{n^k} - \frac{\binom{j}{k}}{\binom{n}{k}} \right] \ \leq\ \frac{2154}{15\,625}\frac{m^2}{n}.
    \end{align}
    with equality for $j=4$ and $m=n=5$.
\end{lemma}
}

\begin{proof}
{
    The summand is maximized when $j=n-1$, thus
    \begin{align}
        \sum_{k=0}^{m} \left[ \frac{j^k}{n^k} - \frac{\binom{j}{k}}{\binom{n}{k}} \right]\ &\leq\ \sum_{k=0}^{m} \left[\frac{(n-1)^k}{n^k} - \frac{\binom{n-1}{k}}{\binom{n}{k}} \right]\\\
        &=\ -(n-1)^{m+1} n^{-m}+\frac{m (m+1)}{2 n}-m+n-1\ \\
        &:=\ S_{m,n}.
    \end{align}
}
{
    Consider
    \begin{align}
        \frac{n}{m^2} S_{m,n}\ &=\ \frac{1}{2}-\frac{2n - 1}{2m} + \frac{(n-1) n^{1-m} \left(n^m-(n-1)^m\right)}{m^2}
    \end{align}
    Treating this as a function in $m$,
    \begin{eqnarray}
        \frac{n^m - (n - 1)^m}{n^{m-1}}\ =\ -\sum_{k=0}^{m-1}(-1)^{m-k}\binom{m}{k}\frac{1}{n^{m-1-k}}\ =\ m + \mathcal{O}(1/n).
    \end{eqnarray}
    As such, up to the leading order we have
    \begin{align}
        \frac{n}{m^2}S_{m,n}\ &=\ \frac{1}{2} - \frac{2n-1}{2m} + \frac{(n-1)m + \mathcal{O}(1)}{m^2}\ =\ \frac{1}{2} - \frac{1}{2m} + \frac{\mathcal{O}(1)}{m^2}.
    \end{align}
    As a function of $m$, $(n/m^2)S_{m,n}$ is therefore an increasing function and maximized at $m=n$. The function $S_{n,n}/n$ achieves a maximum of $2154/15625$ at $n=5$. 
}
\end{proof}

{
Importantly, if $\min_{x\in[0,1]}p(x)$ is small, an arbitrarily large $N\geq m$ is required to have all coefficients positive. As such, a degree 2 polynomial may need to be represented by an arbitrarily large $N$ Bernstein polynomial for the coefficients to be positive. However, it can be converted to a degree $m$ Bernstein polynomial with negative weights and no roots in $[0,1]$.
}

\subsection{Residuals}\label{sec:residuals}
To find the coefficients $w\in\Delta^{M+1}$ and $a\in\mathbb{R}^{N+1}$ in \ref{eq:bernstein_rational_function}, we take the least squares approach and minimize a loss function between the target function and the rational function. However, the strict positivity enforced in the denominator allows for three residual formulations, each with its own benefits.

\textbf{Reweighted Linearized Residuals}: Following the works of Sanathanan and Koerner, one can consider an iterative algorithm where, on the $t$-th iteration, we minimize the reweighted loss
\begin{align}
    \mathcal{L}_{\ell, t}(f, \RR_{N,M}\, |\, a, w)\ =\ \sum_{i=0}^I\, 
    \, \left(\frac{f(x_i)\sum_{m=0}^M\, w_m\, \mathcal{B}_m(x_i)\ -\ \sum_{n=0}^{N}\, a_n\, p_n(x_i)}{\sum_{m=0}^M w_{m}^{t-1}\mathcal{B}_m(x_i)}\right)^2,
\end{align}
where $\{x_i\}_{i=0}^I$ is some partition of $[0, 1]$. Combined with the simplex constraints, this has the benefit of being a quadratic programming problem and thus can be solved by standard solvers.

\textbf{Nonlinear Residuals}: {The strict positivity in the denominator guarantees no poles inside the approximation domain.} Thus, we can consider the nonlinear residuals
\begin{align}
    \mathcal{L}_{r}(f, \RR_{N, M}\, |\, a, w)\ &=\ \sum_{i=0}^I\, \left(f(x_i)\ -\ \frac{\sum_{n=0}^{N}\, a_n\, p_n(x_i)}{\sum_{m=0}^M\, w_m\, \mathcal{B}_m(x_i)}\right)^2. 
\end{align}
While this problem is nonlinear, this formulation represents an accurate formulation of the original approximation problem. To solve this optimization problem, we introduce an iterative scheme in Section \ref{sec:optimization} that allows us to minimize both representations of the residuals.

\subsection{Sobolev-Jacobi Smoothing} {In the presence of noisy data, polynomial regression is known to overfit \cite{intro_to_stat_learning}. One way to counteract this is to add a penalty on the $n$-th derivative of the polynomial. When finding the coefficients, the regularization term has the additional benefit of guaranteeing the corresponding least-squares problem to have a unique solution \cite{Hoerl2000, vanwieringen2023}.} 

To introduce smoothing in rational functions, we consider a penalty on the numerator\footnote{The denominator is not penalized as we find that the simplex constraint is strong enough to maintain smoothness.} induced by the Sobolev inner product
\begin{eqnarray}
    \langle f, g\rangle_\mathcal{H}\ =\ \sum_{k\geq0}\, \lambda_k\, \int_\Omega \, f^{(k)}(x)\, g^{(k)}(x)\, d\rho_k,
\end{eqnarray}
for functions $f$ and $g$, constants $\lambda_k \geq 0$, $\Omega \subseteq \mathbb{R}$ and $\rho_k$ are measures over $\Omega$.

Remarkably, classical orthogonal polynomials are orthogonal with respect to the Sobolev inner product. In particular, the shifted Jacobi polynomials, $\widetilde{P}_n^{(a,b)}(x)$, defined on $[0,1]$ are orthogonal on the domain $[0, 1]$. Taking them as the basis functions for $\PP_N$, with coefficients $\{a_n\}_{n=0}^N$, gives the asymptotic expression
\begin{eqnarray}
    \langle \PP_N, \PP_N\rangle_{\mathcal{H}}\ =\ \sum_{n=0}^N\, a_n^2\, \langle\widetilde{P}_n^{(a,b)}(x), \widetilde{P}_n^{(a,b)}(x)\rangle_\mathcal{H}\ \sim\ \sum_{n=0}^N\, a_n^2 \, \lambda_n\, (2n)^{2n}.
\end{eqnarray}

Thus, in local coordinates (parameterized by Jacobi polynomials, and thus also Gegenbauer, Legendre, and Chebyshev polynomials \cite{orthogonal_polynomials_textbook}), the norm induced by the Sobelev inner product becomes a weighted $\ell^2$ norm on the spectral coefficients, with the coefficients being exponentially damped. That is, an exponential penalty on its spectral coefficients can create a smooth polynomial approximation. To reduce the likelihood of overflow errors in the computations, we take $\lambda_n = 2^{-2n}\, n^{-n}$. This yields what we call the Sobolov-Jacobi smoothing penalty for polynomial approximation
\begin{eqnarray}
    \mathcal{R}(\PP_N)\ =\ \sum_{n=0}^N\, a_n^2\, \, n^{n}.
\end{eqnarray}

\subsection{Multivariate}
So far, we have only considered functions of a single variable. However, the ability to model functions of multiple variables has wide applications in various areas of social sciences (see \cite{spatial_regression_modeling} and the references within). To this end, we generalize our rational function by expressing the numerator and denominator as a tensor product of the respective basis,
\begin{align}\label{eq:multivariate_rational_polynomil}
    \RR(x_1,\ldots, x_s)\ &=\ \frac{\sum_{n_1=0}^{N_1}\cdots\sum_{n_s=0}^{N_s}\, a_{n_1,\ldots, n_s}\, p_{n_1}(x_1)\cdots \, p_{n_s}(x_s)}{\sum_{m_1=0}^{M_1}\cdots \sum_{m_s=0}^{M_s} w_{m_1,\ldots, m_s}\, \mathcal{B}_{m_1}(x_1)\cdots \, \mathcal{B}_{m_s}(x_s)}.
\end{align}
Positivity and normalization of the denominator are enforced by the constraints
\begin{align}
    w_{m_1,\ldots, m_s}\ \geq\ 0\quad\text{and}\quad \sum_{m_1,\ldots, m_s} w_{m_1,\ldots, m_s}\ = \ 1.
\end{align}

Taking the shifted Jacobi polynomials as the basis for the numerator, we generalize the smoothing penalty as
\begin{align}\label{eq:multivariate_smoothing_penalty}
    \mathcal{R}(\NN)\ =\ \sum_{n_1=0}^{N_1}\, \cdots\, \sum_{n_s=0}^{N_s}\ a_{n_1,\ldots, n_s}^2\, n_1^{n_1}\ \cdots\ n_s^{n_s}.
\end{align}

The linearized, reweighted, and nonlinear residuals are easily generalized.

\subsection{Optimization}\label{sec:optimization}
To find the coefficients for the rational Bernstein Denominator approximation, $\RR(x_1,\ldots, x_s)$, we consider the constrained optimization problem
\begin{eqnarray}
    \min_{{a\in\mathbb{R}^{{N}'},\ w\in\Delta^{{M}'}}}\ \mathcal{L}(f, \RR\, |\, a, w)\ +\ \lambda\, \mathcal{R}(\NN),
\end{eqnarray}
where ${N}'=s + N_1+\cdots+N_s$, ${M}'=s + M_1+\cdots+M_s$, some smoothing strength $\lambda \geq 0$, and $\mathcal{L}$ can be the reweighted, or the nonlinear residuals.

The reweighted residuals have the advantage of being a standard quadratic programming problem, even with the regularization term. However, this formulation often results in worse fits than minimizing the nonlinear residuals directly. 

We propose the following iterative scheme to solve for the nonlinear case (and can be applied to the linearized and reweighted cases). Let $(a^t, w^t)$ be the value of $(a, w)$ at iteration $t$. For the next iteration, we use the scheme proposed by Chok and Vasil \cite{james_vasil}
\begin{align}
    & w^{t+1}_i\ =\ w^t_i\, -\, \eta_t d^t_i\quad \text{with}\\
    \notag&d^t_i\ =\ w^t_i(\nabla_w\mathcal{L}^t_i\, -\, w^t\cdot\nabla_w \mathcal{L}^t),\quad\text{for}\quad 0 < \eta_t \leq \eta_{t,\max},\\
    \notag& \eta_{t,\max}^{-1}\ =\ \max_i(\nabla_{w}\mathcal{L}^t_i\, -\, w^t\cdot \nabla_w\mathcal{L}^t),
\end{align}
and $\nabla_w\mathcal{L}^t = \nabla_w\mathcal{L}(f, \RR\, |\, a^t, w)$. This iteration scheme, as shown in their paper, maintains that $w^{t+1}\in\Delta^{M'}$ if $w^t\in\Delta^{M'}$. Then $a^{t+1}$ is chosen as the solution to
\begin{eqnarray}
    a^{t+1} = \underset{a\in\mathbb{R}^{N'}}{\arg \min}\ \mathcal{L}(f, \RR\, |\, a, w^{t+1})\ +\ \lambda\, \mathcal{R}(\NN),
\end{eqnarray}
which is a standard quadratic programming problem that has an explicit solution. In the first iteration, $w^{0}=\frac{1}{{M}'}\mathbbm{1}$ and $\mathbbm{1}\in\mathbb{R}^{{M}'}$ is a vector filled with ones. 

\begin{remark}
    The iteration for $w^{t+1}$ can be performed once or multiple times. However, we find that running it once yields better results.
\end{remark}

\subsubsection{Hot-Start}
Iterative methods for nonlinear problems, \textit{i.e.} the nonlinear residuals, are often sensitive to initial conditions. To speed up the iterative procedure and potentially converge to a better solution, we propose to \textit{hot-start} the iterative algorithm with a solution obtained from an alternative rational approximation, typically the SK or AAA algorithm, or the rational Bernstein denominator algorithm with reweighted linearized residuals.

We first fit the alternative rational approximation and the rational Bernstein denominator with reweighted linearized residuals. The denominators are converted to a Bernstein basis and the coefficients are projected onto the simplex. Once projected, they may contain poles at $x=0$ or $x=1$. If both do not contain poles, then the one with the best error is used to hot-start the iteration. Otherwise, if only one does not contain poles, that one is used to begin the iteration. If both contain poles, then we take $\frac{1}{M'}\mathbbm{1}$ as the initial condition.

\section{Numerical Results}
{All our code is in Python and can be found on our Github\footnote{https://github.com/infamoussoap/RationalFunctionApproximation}. Our code implements the Bernstein Denominator algorithm and other rational function approximation methods, like the AAA and SK algorithms.}

\subsection{Scaling of Bernstein polynomials} 
{As seen in Theorem \ref{theorem:polynomial_in_bernstein}, a strictly positive polynomial $p(x)$ on $[0,1]$ can be represented by a degree $N$ Bernstein polynomial with strictly positive coefficients. However, $N$ scales like $\mathcal{O}(\frac{1}{\min_{x\in[0,1]} p(x)})$. }

\begin{figure}
    \centering
    \includegraphics[width=0.98\linewidth]{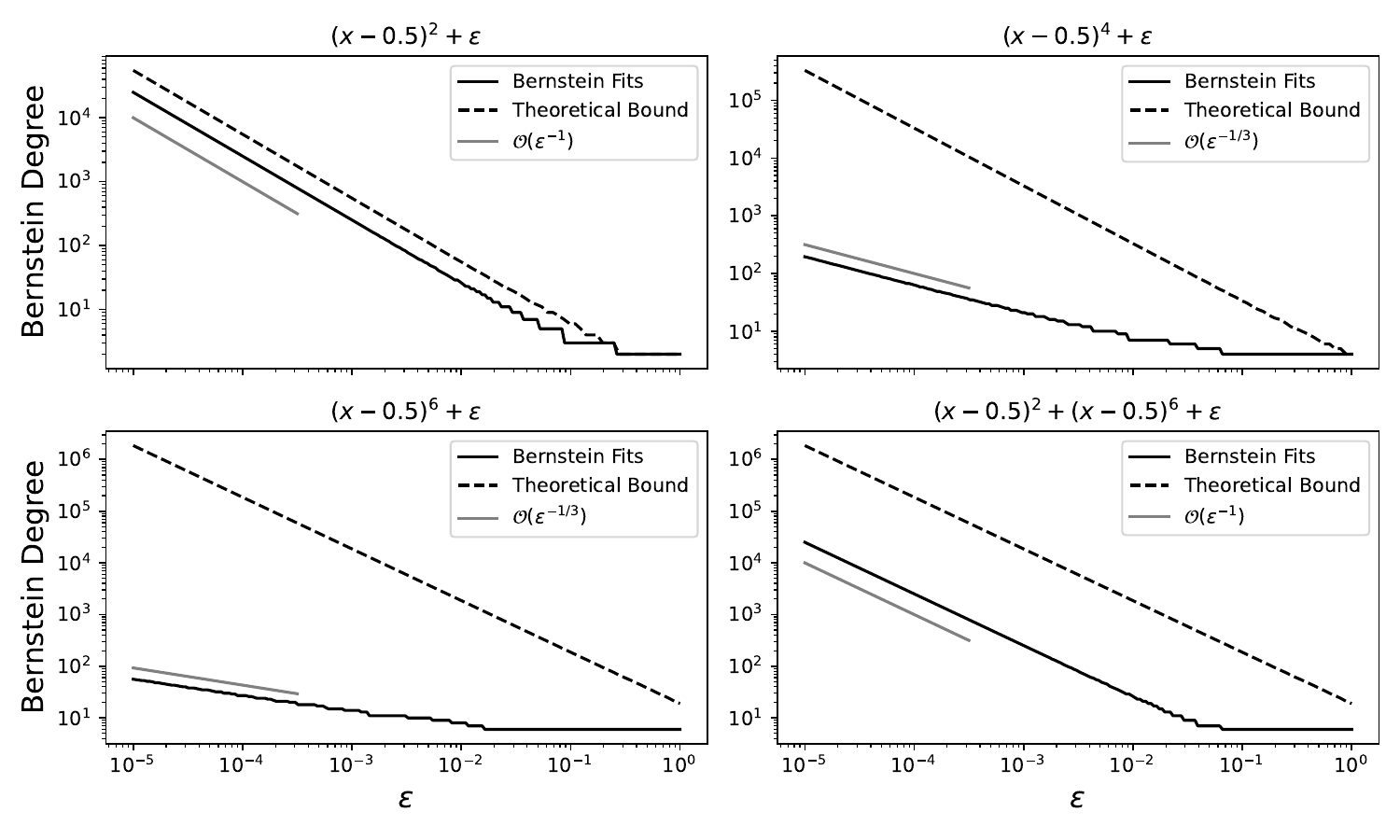}
    \caption{Minimum Bernstein degree to have positive coefficients when representing various polynomials. This is compared to the theoretical bound (\ref{eq:positive_bernstein_bound}).}
    \label{fig:positive_bernstein_convergence}
\end{figure}

{
To that end, we consider the polynomials $(x-0.5)^k+\varepsilon$ for $k=2,4,6$, and $(x-0.5)^2 + (x-0.5)^6+\varepsilon$ for various values of $\varepsilon$. Using the conversion formula (\ref{eq:monomials_to_bernstein_coef}), we find the minimum degree $N$ Bernstein polynomial such that all coefficients $>10^{-10}$. This is compared to the theoretical bound (\ref{eq:positive_bernstein_bound}) and can be seen in Fig. \ref{fig:positive_bernstein_convergence}. }

{
The minimum degree $N$ lies below the theoretical bound and is sharp for $(x-0.5)^2+\varepsilon$ when $0.1\leq\varepsilon\leq1$. For $(x-0.5)^k+\varepsilon$, when $k=4,6$, $N$ scales like $\varepsilon^{-1/2k}$. However, this is a special case, as we see the expected $\mathcal{O}(1/\varepsilon)$ rate for $(x-0.5)^2+(x-0.5)^6+\varepsilon$.}

\subsection{Differential Equations}\label{sec:differential_equations}
As outlined in the introduction, consider solving a general non-constant coefficient differential equation
\begin{align}\label{eq:general_ode}
    \sum_{k=0}^K c_k(x) D^k u(x)\ =\ f(x),
\end{align}
where $f(x)\in \mathcal{P}(n)$, for some positive integer $n$. Spectral methods numerically solve these differential equations by approximating the solution as a polynomial
\begin{align}
    u\ \approx\ u_N\ =\ \sum_{i=0}^N\, a_i\, \phi_i(x),
\end{align}
for some orthogonal basis functions $\{\phi_i\}$ known as trial functions, with $\{a_i\}$ to be found. Appropriate basis functions convert the spatial derivatives into sparse matrices \cite{dedalus, spectral_methods_sparse}, turning (\ref{eq:general_ode}) into a linear system of equations that can be solved efficiently.

When $p(x)$ and $q(x)$ are not constant, numerical spectral solvers, like Dedalus \cite{dedalus}, approximate them in terms of orthogonal basis functions $\{\psi_i\}$, known as test functions. This converts them into banded matrices approximately the size of the number of basis functions used. Once multiplied with the derivative matrices, one may no longer have a sparse system. Thus, one must strike a balance between accurate approximation while using as few coefficients as possible.

A natural proposal would be to use rational approximations instead. Concretely, for (\ref{eq:general_ode}), approximate the non-constant coefficients $c_k(x)$ as rational functions
\begin{align}
    c_k(x)\ \approx\ \frac{p_{k, N}(x)}{q_M(x)},
\end{align}
where $p_{k, N} \in \mathcal{P}(N)$ and $q_M \in \mathcal{P}(M)$. If $q_M(x)\neq 0$ in the approximation interval, we can instead solve
\begin{equation}
\sum_{k=0}^{K} p_{k,N}(x) \, D^{k} u(x) \ = \ q_{M}(x) \, f(x). \label{common denominator}
\end{equation} 
In this situation, one wants a guarantee of having a strictly positive denominator. Otherwise, the numerical spectral solver will return the wrong results. As such, the strict positivity guaranteed by the Bernstein Denominator algorithm formulation provides a crucial building block for this problem.

This gives two different situations: i) a lower degree of sparsity in the rational approximation while maintaining similar accuracy as the polynomial approximation or ii) a similar degree of sparsity in the rational approximation while having higher accuracy than the polynomial approximation. 

\subsubsection{{Exponential Approximation}}

{In this section, we consider the non-constant-coefficient wave model,
\begin{align}
\label{eq:bessel_single_non-constant_coefficient}
    (\lambda\, e^{2az}\ -\ m^2)\, {y}(z)\ +\ \frac{1}{a^2}{y}''(z)\ =\ 0, \qquad \text{for} \qquad z \in  [0,1].
\end{align}
The system (\ref{eq:bessel_single_non-constant_coefficient}) happens to have solutions in terms of exponentially re-parameterized Bessel functions
\begin{align}
y(z) = c_{1}\, J_{m}(\sqrt{\lambda} e^{az}) + c_{2}\, Y_{m}(\sqrt{\lambda} e^{az}), \label{eq:bessel_equation_general_sol}
\end{align}
where $J_m(x)$ and $Y_m(x)$ are Bessel functions of the first and second kind, respectively, and $c_1,c_2\in\mathbb{R}$.
}

{This example comes from a wave-physics problem modeling the interior of stars and planets \cite{bessel_exponential_paramterization}. We choose this simple example because we want to satisfy three basic requirements: 
\begin{itemize}
    \item We want the ultimate solution in analytical form to compute errors in with numerical approximation scheme. 
    \item We want the model problem to have all the challenges of a typical state-of-the-art full-production simulation, but without the many expensive complications that come with three-dimensional high-performance computations. 
    \item Specifically, we want to single out the dependence on large-dynamic-range non-constant coefficients.  
\end{itemize}}

{For the first point, the analytical solution is directly apparent by design. However, we note that the original authors of the simple wave model did not construct the problem to be analytically solvable. The original model was derived to include semi-realistic sharp parameter transitions in the wave-medium properties; the Bessel solution was serendipitous.}

{For the second point, we note that a problem like (\ref{eq:bessel_single_non-constant_coefficient}) occurs naturally as a plausible model for a single mode out of a large number of decoupled (or semi-decoupled) linear boundary-value problems, each forced independently by time-lagged nonlinear terms \cite{dedalus}. That is, (\ref{eq:bessel_single_non-constant_coefficient}) assumes the general structure of a single part of modern time-integrated wave-convection physics problems in massive stars \cite{Anders2023} or a decoupled eigenmode problem in solar-interior physics \cite{Vasil2024}.}

{To the third point, we we set $a=8$ and $m=2$, producing the the non-constant coefficient $e^{16 z}$, with an overall dynamic range across the domain of $e^{16} \approx 9\times 10^{6}$. The exponential profile is far from a simple polynomial. Moreover, by its physical nature, the coefficient must remain everywhere \textit{positive}. A spectral polynomial approximation will typically achieve high accuracy by creating errors with uniform oscillation. Given the large dynamic range those oscillating error can make the approximation become negative where the values of the original function are small. While our specific rational approximation only enforces positivity of the denominator, the denominator acts as an overall weight that mollifies the total dynamic range, making the uniform oscillations from the numerator less likely to cause negative values.}

{\textbf{Experiment Details}: For increasing values of $n$, we compare the approximations of a $\mathcal{P}(2n)$ with $\mathcal{R}(n,n)$ AAA and rational Bernstein denominator. All methods were fitted on a uniform grid on $[0, 1]$ with $4096$ points. We record the maximum absolute error (MaxAE) and check if poles exist in $[0, 1]$. Each method was run five times to record the average time taken\footnote{We run all experiments in Section \ref{sec:differential_equations} on a 2021 MacBook Pro with an Apple M1 Pro chip.}.}

{The $\mathcal{R}(n,n)$ rational Bernstein denominator is initialized with a $\mathcal{R}(m, m)$ AAA where $m=\min(11,n)$ as we find that $n>11$ causes the AAA algorithm to produce poles.\footnote{For the exponential function, we find that the AAA produces better initial conditions than the SK algorithm.} We also minimize the nonlinear residuals.}

\begin{figure}
    \centering
    \includegraphics[width=0.98\linewidth]{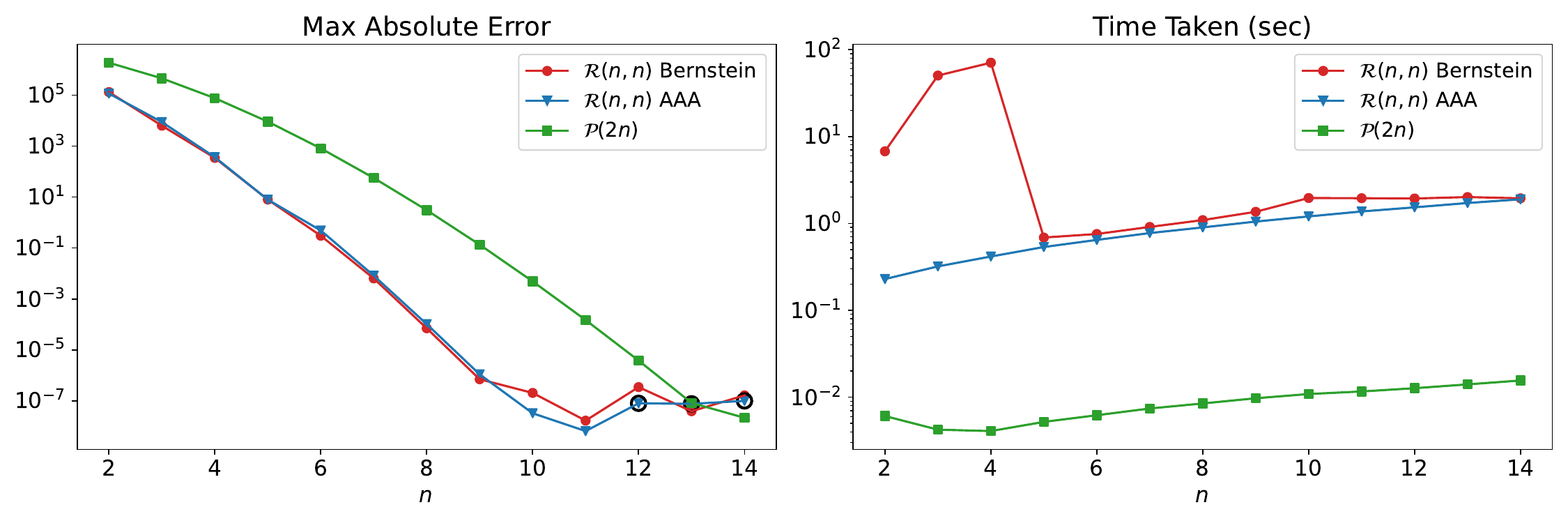}
    \caption{Maximum absolute error between the true function and the fitted methods on $\exp(16 x)$. The black circle indicates fits with poles inside the approximation interval $[0, 1]$.}\label{fig:exp_convergence_plots}
\end{figure}

{The results can be seen in Fig. \ref{fig:exp_convergence_plots}.}

{\textbf{Results}: In accuracy, the AAA and rational Bernstein denominator outperform the equivalent polynomial approximation. The AAA and rational Bernstein denominator also share similar levels of accuracy for $n\leq 9$, after which the AAA algorithm yields better accuracy. The polynomial approximation reaches a minimum of $3.82\times10^{-6}$ at $\mathcal{P}(24)$, while the AAA and rational Bernstein denominator, respectively, reach a minimum of $8.01\times10^{-9}$ and $1.68\times10^{-8}$ at $\mathcal{R}(11,11)$. After this, the AAA algorithm produces a pole for $n\geq12$ while the rational Bernstein denominator does not produce poles. However, the rational Bernstein denominator increases in error to $3.46\times10^{-7}$.}

{In time, the polynomial approximation is solved faster than the AAA which, in turn, is faster than the rational Bernstein denominator. However, when $n=2,3,4$, the rational Bernstein denominator requires 10-100 times more time than the equivalent AAA with little benefit to accuracy. This is due to the rational Bernstein denominator minimizing the $\ell^2$ norm while we report the $\ell^\infty$ norm. Moreover, the AAA yields an initial condition `far' from a minimum. As such, more time is spent in the iterative solver to find a better solution. When $n\geq5$, the AAA yields an initial condition near a minimum for the rational Berenstein denominator. This results in less time spent in the iterative solver.}

\subsubsection{Numerical Solution}
\begin{figure}
    \centering
    \includegraphics[width=0.98\linewidth]{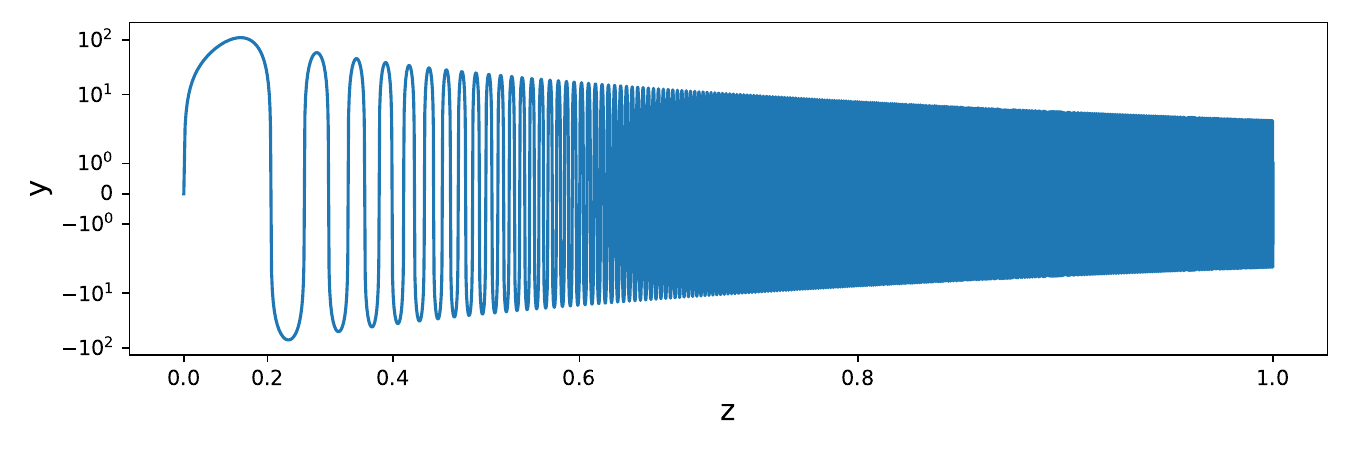}
    \caption{Solution to Bessel's differential equation with $m=2$, $a=8$ and initial conditions ${y}(0)=0$, ${y}(1) = 1$.}
    \label{fig:bessel_true_solution}
\end{figure}

{We now consider solving the exponential parameterized Bessel differential equation with boundary conditions $y(0)=0$ and $y(1)=1$. The exact (highly oscillatory) solution can be found by solving for $c_1$ and $c_2$ in (\ref{eq:bessel_equation_general_sol}). This can be seen in Fig. \ref{fig:bessel_true_solution}.}

{\textbf{Experimental Details}: For a given $n$, the non-constant coefficient is approximated with a $\mathcal{P}(2n)$ polynomial, $\mathcal{R}(n,n)$ AAA (with cleanup on) and rational Bernstein denominator, as performed in the previous section.\footnote{The AAA algorithm is fitted using the Barycentric form and then converted to a polynomial.} These approximations are then fed into Dedalus \cite{dedalus} using a Chebyshev-$\tau$ scheme with $2^{14}=16\,384$ spectral modes. The $\ell^2$ norm between the true and numerical solution is recorded.} 

{In a full state-of-the-art production calculation, an equation like (\ref{eq:bessel_single_non-constant_coefficient}) would comprise the left-hand side of a time-stepping scheme. To this end, the resulting matrix would be factorized initially and solved repeatedly for different right-hand-side over (e.g.) $\approx 10^{6}$ time steps.}

{The error in $\mathcal{P}(2n)$ is roughly comparable to $\mathcal{R}(n,n)$. However, the non-constant coefficient off-diagonal \textit{matrix band-width} of $\mathcal{P}(2n)$ is exactly double that of $\mathcal{R}(n,n)$. After computing a banded matrix factorization, the forward-backward solve complexity scales as $\mathcal{O}(b^{3})$, where the bandwidth, $b = n$ for $\mathcal{P}(2n)$ and $b=n$ for $\mathcal{R}(n,n)$. \textit{This is the cost we are concerned with mitigating.}}

{Moreover, one might be tempted to try a rational approximation with different degrees, $\mathcal{R}(n,m)$ where $n\ne m$. While this might be useful for isolated approximations, it defeats the purpose for the goal of minimizing matrix band-width. In (\ref{common denominator}) we need to multiply through by the denominator. The goal is to balance the band-width right- and left- hand sides while minimizing approximation error. If we used smaller denominator degree, we would leave potential decrease in the numerator degree on the table. Alternatively, increasing the  denominator degree beyond the numerator degree simply reverses the issue.}

{The Lower-Upper (LU) matrix factorization is the most obvious algorithm. However, in practice, the LU factorization produces dense fill-in for $\tau$-corrected spectral methods like we consider here. Alternatively, we find excellent sparsity preservation with an Upper-Lower (UL) matrix factorization, which we compute via an LU factorization on the transpose of the matrix. We use the same SuperLU \cite{SuperLU} library that Dedalus calls for production calculations. In general, the cost to perform the factorization is greater than the solve time after obtaining $U$ and $L$ (which is the point). In practice, the factorization costs are amortized over many more individual solves against multiple right-hand sides in a time-stepping computation.  
We run the UL factorization and corresponding solve 100 times to record the average time taken.}

\begin{figure}
    \centering
    \includegraphics[width=0.98\linewidth]{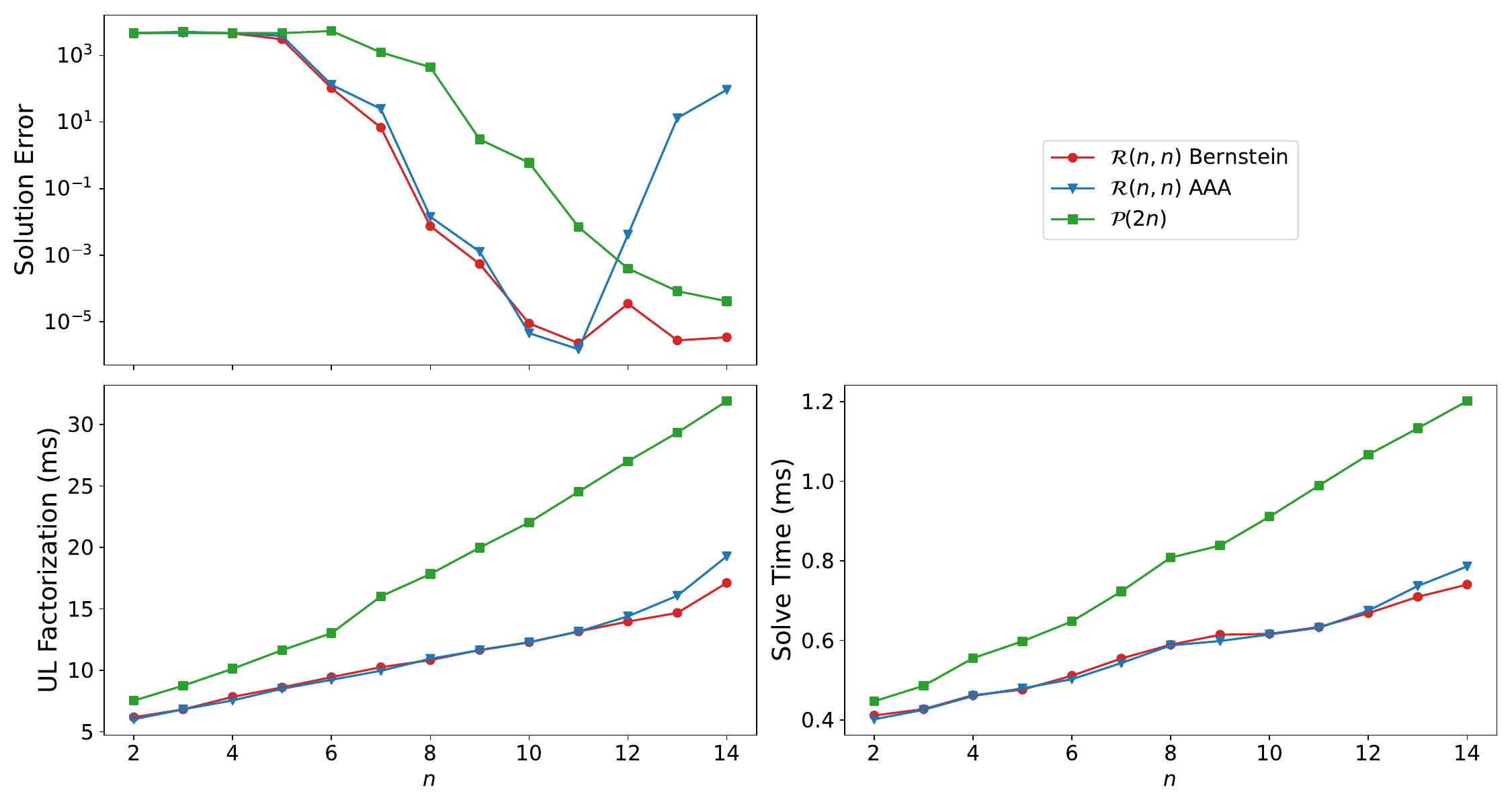}
    \caption{Solution errors and timings when solving the non-constant-coefficient wave model using a Chebyshev-$\tau$ scheme with polynomial and rational approximations of the non-constant coefficient.}
    \label{fig:lu_error_solve_times}
\end{figure}

{The results can be seen in Fig. \ref{fig:lu_error_solve_times}. We also record the sparsity patterns in the matrix before and after the UL factorization is performed in Table \ref{tab:sparsity_patterns}.}

{\textbf{Results}:  
For $n=2,3,4$ and $5$, all three approximation methods yield similar errors. There is a noticeable divergence for $n\geq 6$ where the rational approximations far outperform the polynomial counterparts. In particular, the polynomial reaches a minimum error of $4.14\times10^{-5}$ at $\mathcal{P}(28)$ while the AAA and rational Bernstein denominator obtain a minimum error of $1.50\times10^{-6}$ and $2.28\times10^{-6}$, respectively, at $\mathcal{R}(11, 11)$.}

{For $n\geq 12$, the AAA algorithm experiences a notable increase in error, producing solutions worse than the corresponding polynomial approximation. However, this is to be expected as the AAA algorithm produces poles in $[0, 1]$, highlighting the importance of pole-free rational approximations in $[0, 1]$.}

{In time, we see that the polynomial approximation always takes longer than the corresponding rational approximations for both the UL factorization and solve times. In particular, the runtime scales with the matrix bandwidth size and not the degree of freedom of the non-constant coefficient approximation, as expected. That is, a $\mathcal{R}(n,n)$ approximation and a $\mathcal{P}(n)$ or $\mathcal{P}(n+1)$ approximation take similar lengths of time. This fact is corroborated by the sparsity patterns of the original matrix and the corresponding UL decomposition, where the $\mathcal{R}(n, n)$ rational approximations have a similar count of non-zero entries as a $\mathcal{P}(n)$ or $\mathcal{P}(n+1)$ approximations.}

\begin{table}[]
    \centering
    \caption{Number of non-zero elements when non-constant-coefficient wave model is discretized using a Chebyshev-$\tau$ scheme equation and its UL decomposition. In brackets are the inflation factors, the non-zero elements in the UL factorization divided by the non-zero elements in the original matrix.}
    \label{tab:sparsity_patterns}
    \begin{tabular}{l|ccc|ccc}
\toprule
 & \multicolumn{3}{c|}{Orig. Matrix} & \multicolumn{3}{c}{UL Factorization} \\
 \midrule
 &$\mathcal{P}(2n)$ & $\mathcal{R}(n, n)$ &$\mathcal{R}(n, n)$ &$\mathcal{P}(2n)$ & $\mathcal{R}(n, n)$ &$\mathcal{R}(n, n)$ \\
 $n$& Polynomial & AAA & Bernstein & Polynomial & AAA & Bernstein \\
\midrule
2 & 245,717 & 180,209 & 180,209 & 282,350 (1.15) & 219,290 (1.22) & 221,181 (1.23) \\
3 & 311,223 & 212,969 & 212,969 & 342,637 (1.10) & 255,245 (1.20) & 256,088 (1.20) \\
4 & 376,721 & 245,727 & 245,727 & 408,233 (1.08) & 294,002 (1.20) & 295,377 (1.20) \\
5 & 442,211 & 278,483 & 278,483 & 472,384 (1.07) & 331,933 (1.19) & 331,931 (1.19) \\
6 & 507,693 & 311,237 & 311,237 & 537,388 (1.06) & 369,919 (1.19) & 374,354 (1.20) \\
7 & 573,167 & 343,989 & 343,989 & 602,667 (1.05) & 418,370 (1.22) & 418,324 (1.22) \\
8 & 638,633 & 376,739 & 376,739 & 668,136 (1.05) & 435,827 (1.16) & 436,648 (1.16) \\
9 & 704,091 & 409,486 & 409,487 & 733,587 (1.04) & 473,993 (1.16) & 473,977 (1.16) \\
10 & 769,541 & 442,231 & 442,233 & 799,033 (1.04) & 496,629 (1.12) & 496,628 (1.12) \\
11 & 834,983 & 474,975 & 474,977 & 864,471 (1.04) & 530,460 (1.12) & 530,808 (1.12) \\
12 & 900,417 & 507,716 & 507,717 & 929,901 (1.03) & 599,620 (1.18) & 563,204 (1.11) \\
13 & 965,843 & 540,449 & 540,455 & 995,323 (1.03) & 648,350 (1.20) & 596,279 (1.10) \\
14 & 1,031,261 & 573,181 & 573,191 & 1,060,737 (1.03) & 697,364 (1.22) & 629,011 (1.10) \\
\bottomrule
\end{tabular}
\end{table}

{
\subsection{Barycentric vs Bernstein}
We now consider the numerical convergence of AAA\footnote{While our experiments use our Python implementation of AAA, the results were cross-checked with the original implementation in Matlab.} and the Floater and Hormann poleless Barycentric form\footnote{To ensure our implementation of the poleless Barycentric form was correct, our results were cross-checked with the results in their paper \cite{barycentric_poleless}.} and Bernstein Denominator algorithm for the functions $|2(x - 0.5)|$, $\exp(-x)\sin(16x^2)$, $1/[(x-0.01)(x+1.01)(x^2+0.01)]$, and $1/(1 + [10(x - 0.5)]^2)$ in the domain $[0, 1]$. }

{
\textbf{Experimental Details}: For increasing values of $n$, we compare the approximations of a $\mathcal{R}(n, n)$ AAA (with clean-up on), $\mathcal{R}(n + 1, n + 1)$ poleless Barycentric with $d=3$, $\mathcal{R}(n, n)$ rational Bernstein Denominator, and a $2n$-degree polynomial. All methods, except the poleless Barycentric, were fitted on $\{(i/1000, f(i/1000))\}_{0\leq i\leq 1000}$. The poleless Barycentric was fitted on $\{(i/n, f(i/n))\}_{0\leq i\leq n}$, that is $n + 1$ equally spaced points\footnote{Like AAA, one can choose the location of $x_k$ iteratively by selecting locations with the biggest error to the true function value. However, we find that for poleless barycentric forms, this performs worse than choosing points equally spaced} on $[0, 1]$.}

{
For all methods, the maximum absolute error (MaxAE) between the approximated function and the target function $f(x)$ is assessed at the locations $\{i/5000\}_{i=0}^{5000}$. Each fit was checked to see if it contained poles within the interval $[0, 1]$. }

{
For the Bernstein Denominator algorithm, we do not apply Sobolov-Smoothing, but apply hot-start and solve for the nonlinear residuals.}

{
All fits were run on an Intel Xeon Gold 6248 CPU.}

\begin{figure}[ht]
    \centering
    \includegraphics[width=0.99\linewidth]{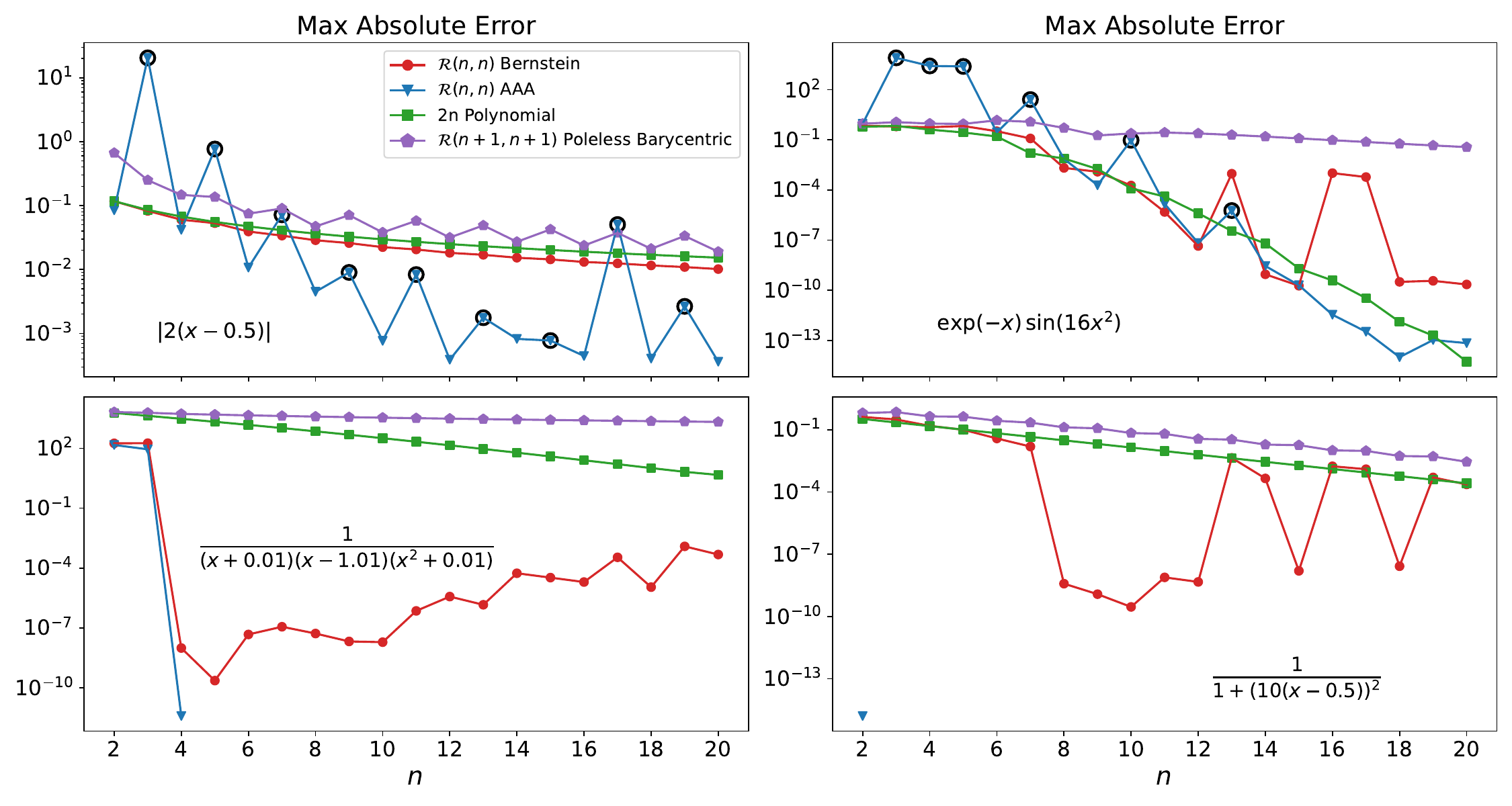}
    \caption{Maximum absolute error between the true function and the fitted methods on various functions. The black circle indicates fits with poles inside the approximation interval $[0,1]$. For $1/(1 + [10(x-0.5)]^2)$, AAA has converged at $n=2$.}
    \label{fig:aaa_bernstein_convergence_plots}
\end{figure}

\begin{figure}[ht]
    \centering
    \includegraphics[width=0.99\linewidth]{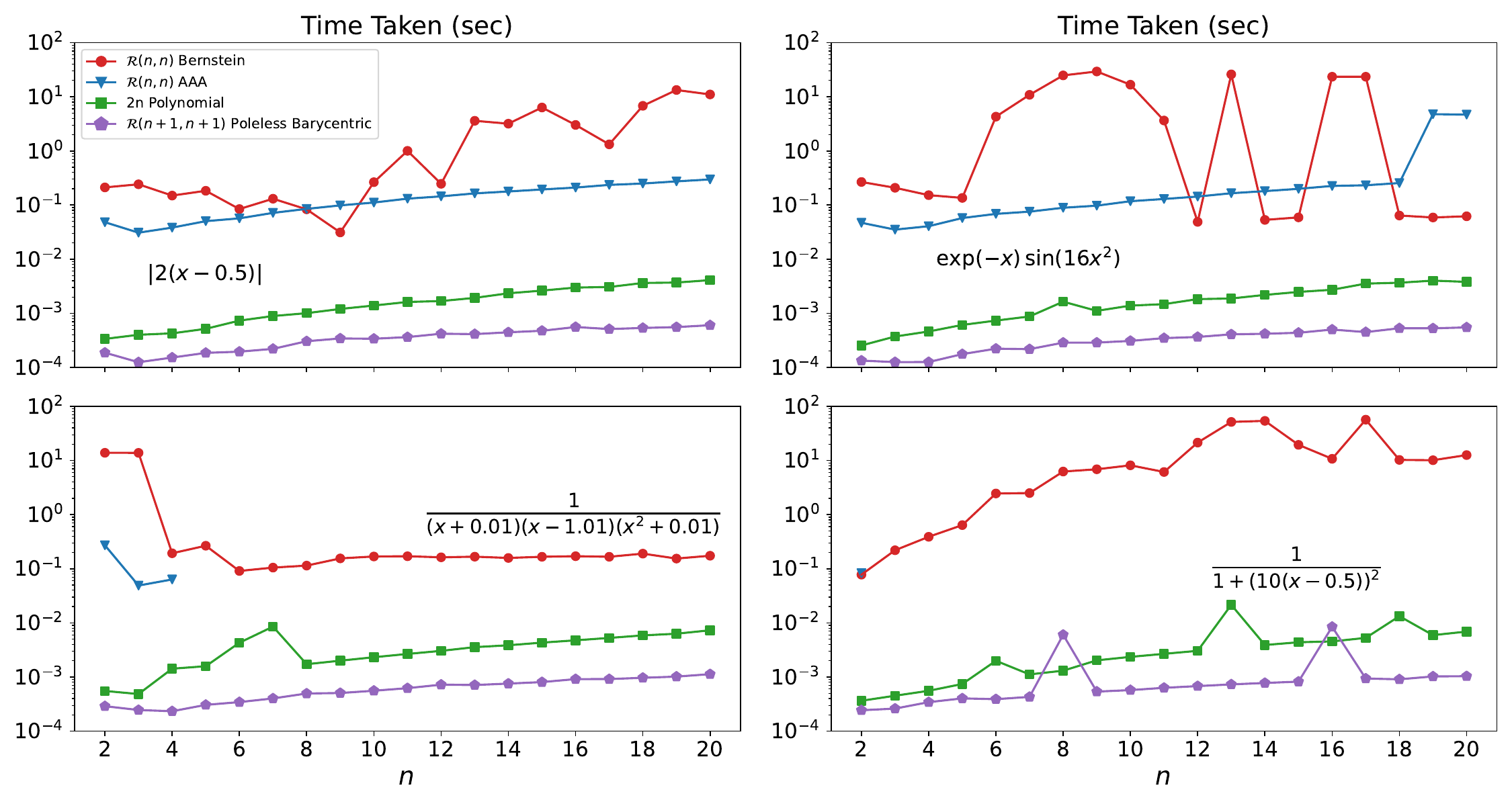}
    \caption{Time taken when fitting $2n$ Polynomial, $\mathcal{R}(n, n)$ AAA (with cleanup on), $\mathcal{R}(n + 1, n + 1)$ Poleless Barycentric, and $\mathcal{R}(n, n)$ rational Bernstein Denominator fits on various functions.}
    \label{fig:aaa_bernstein_time_plots}
\end{figure}

{The numerical convergence plots can be seen in Figure \ref{fig:aaa_bernstein_convergence_plots}, and Figure \ref{fig:aaa_bernstein_time_plots} displays the run time of the algorithms.}

{\textbf{Results}: The MaxAE of the $\mathcal{R}(n+1,n+1)$ Poleless Barycentric form is consistently outperformed by a $2n$ degree polynomial. However, as shown in their paper \cite{barycentric_poleless}, we note that it performs comparably to cubic splines.}

{The Bernstein Denominator consistently outperforms the equivalent Poless Barycentric form. However, its approximations often fluctuate after $n=10$. This is likely due to converging to a non-global minimum, showing the difficulty of solving the nonlinear residuals. We also see that the Bernstein Denominator algorithm generally converges at the same rate, and sometimes better, than the equivalent $\mathcal{P}(2n)$ degree polynomial. Thus, the expected factor of 2 improvement is attained.}

{As expected, the AAA algorithm consistently outperforms all methods. However, as noted in the original AAA paper \cite{AAA}, it produces poles for $|2(x-0.5)|$ for odd $n$. Similarly, for $\exp(-x)\, \sin(16x^2)$, the AAA algorithm produces poles for $n = 3, 4, 5, 7, 10$ and $13$, but doesn't produce poles when $n\geq 14$. Moreover, the AAA and Bernstein Denominator converge at similar rates to a polynomial approximation for $\exp(-x)\sin(16 x^2)$.}

{Since the AAA and Bernstein Denominator algorithms require iterative solvers, it is unsurprising that they consistently take the longest time to run. However, the Bernstein Denominator algorithm consistently takes the longest to run, with runtimes ranging from $0.031$-$56.75$ seconds. Interestingly, a long run time does not necessarily correspond to more accurate solutions. This is particularly prominent in $e^{-x}\sin(16x^2)$ in which $n=16,17$ have a runtime of $24.5$ seconds but yield worse MaxAE than for $n=18,19,20$ which took around $0.06$ seconds. This is likely due to the hot start converging near the optimal value for the latter.}

\subsection{Smoothing Splines vs Bernstein}
Smoothing splines are often a popular choice when modeling spatial data \cite{spline_tensor_2, spline_tensor_3, wood_spline_tensor}, as they provide the flexibility of high-degree polynomials without spurious oscillations. As such, we compare the Bernstein Denominator algorithm to smoothing splines in functions of two variables. 

{\textbf{Experimental Details}: We consider the following functions
\begin{eqnarray}
    &&g_1(x, y)\ =\ \frac{1}{\exp(x/2)}\exp(-x^2)\exp\left(-\frac{y^2}{2\exp(x)}\right), \\
    &&g_2(x, y)\ =\ 20\exp(-0.2\sqrt{0.5[x^2 + y^2]}) + \exp(0.5[\cos 2\pi x + \cos 2\pi y]),
\end{eqnarray}
respectively named Neal's funnel \cite{Neal2003} and Ackley's function \cite{Ackley1987}. We approximate the rescaled functions 
$f_1(x, y)=g_1(8(x-0.5), 8(y-0.5))$ and $f_2(x,y) = g_2(2(x-0.5), 2(x-0.5))$ in noiseless and noisy datasets inside the domain $[0,1]\times[0,1]$.}

For the noiseless data, we compare the numerical convergence on the data \\
$\{(j/50, k/50, f_i(j/50, k/50)\}_{0\leq j, k\leq 50}$ for $i=1$ and $2$. To compare similar degrees of freedom, for a given $n$, we take the Bernstein Denominator to be degree $n$ in both $x$ and $z$ for the numerator and denominator, yielding $2(n + 1)^2 - 1$ degrees of freedom. This is compared to a tensor product of penalized splines with $\Big{\lceil}\sqrt{2(n + 1)^2 - 1}\Big{\rceil}$ basis splines in each variable, yielding $\Big{\lceil}\sqrt{2(n + 1)^2 - 1}\Big{\rceil}^2 + 1$ free coefficients, including the constant term. 

\begin{figure}[ht]
    \centering
    \includegraphics[width=0.95\linewidth]{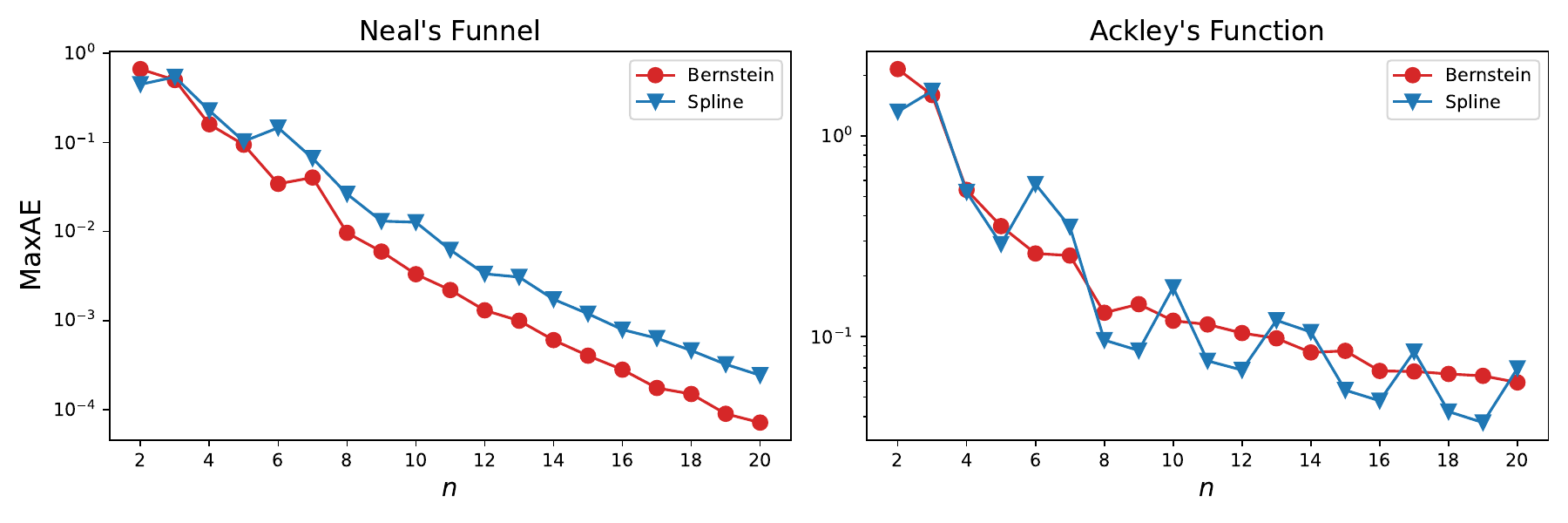}
    \caption{Numerical convergence on 2D functions without noise for a Multivariate rational Bernstein Denominator compared with a tensor product smoothing spline with an equivalent degree of freedom.}
    \label{fig:gam_bernstein_convergence}
\end{figure}

Smoothing and hot-start are not applied for the rational Bernstein denominator and we solve for the nonlinear residuals.
The MaxAE is recorded between the approximation and the true function $f_i$ at the sample points $\{(j/100, k/100)\}_{0\leq j, k\leq 100}$. The results are in Fig. \ref{fig:gam_bernstein_convergence}.

\begin{figure}[ht]
    \centering
    \includegraphics[width=0.98\linewidth]{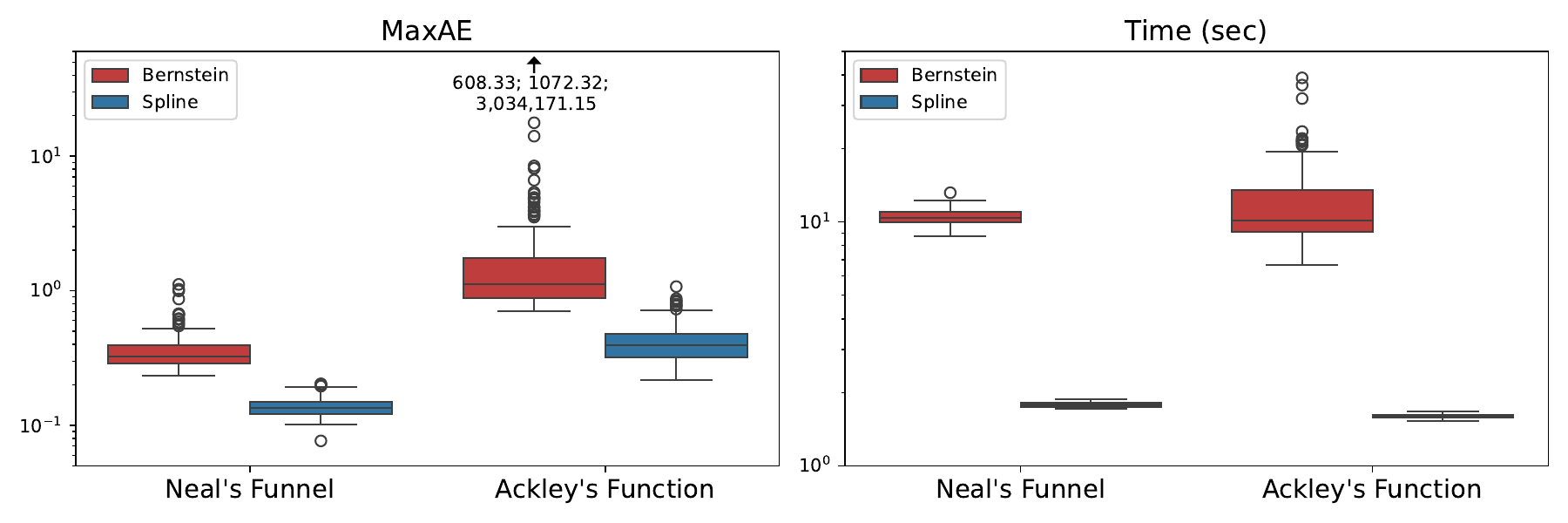}
    \caption{The MaxAE and time taken to fit Sobolov-Smoothed Multivariate rational Bernstein Denominator algorithm and tensor product smoothing splines, with cross-validation, when fitted on noisy data.}
    \label{fig:gam_bernstein}
\end{figure}

For the noisy data, we generate the dataset $\{(x_j, z_j, y_j)\}_{1\leq j\leq 1000}$, where $x_j, z_j \sim U(0, 1)$ and $y_j \sim f_i(x_j, z_j) + \mathcal{N}(0, 0.1^2)$. Two hundred replicate simulates of the dataset were performed. Using cross-validation on each replicate, we selected the smoothing penalty and the degree in the Multivariate rational Bernstein Denominator/degree of freedom for penalized splines. We do not apply hot-start and solve for the nonlinear residuals for the Multivariate rational Bernstein Denominator. We compute the MaxAE between the approximated function and the true function $f_i$ at the sample points $\{(j/100, k/100)\}_{0\leq j, k\leq 100}$. The results are in Fig. \ref{fig:gam_bernstein}.

{\textbf{Results}: In the noise-free case, the Bernstein denominator algorithm converges faster than splines for Neal's funnel, with similar convergence rates for Ackley's function. In time, the Bernstein algorithm takes much longer to produce the approximation. In particular, it takes a minimum of $0.51$ and $0.31$ seconds with a maximum of $130.05$ and $35.32$ minutes for Neal's funnel and Ackley's function respectively. In contrast, splines take a minimum of $0.03$ for both functions and a maximum of $3.09$ and $3.44$ seconds, respectively.}

{In the noisy case, the Sobolov-Smoothed Bernstein denominator algorithm produces fits with a wide variance of errors. This is particularly noticeable for Ackley's function which yields a maximum error of $3.03\times10^{6}$. In contrast, splines have a maximum error of $1.07$. }

{Similarly, the median error of the rational Bernstein denominator fits is higher than the spline fits. The rational Bernstein denominator has a median error of $0.40$ and $1.12$, while the spline has a median error of $0.13$ and $0.32$ for Neal's funnel and Ackely's function respectively. In time, the rational Bernstein denominator takes a median of $10.42$ and $10.16$ seconds, while the spline takes a median of $1.77$ and $1.59$ on Neal's funnel and Ackley's function respectively. }

\begin{figure}[ht]
    \centering
    \includegraphics[width=0.95\linewidth]{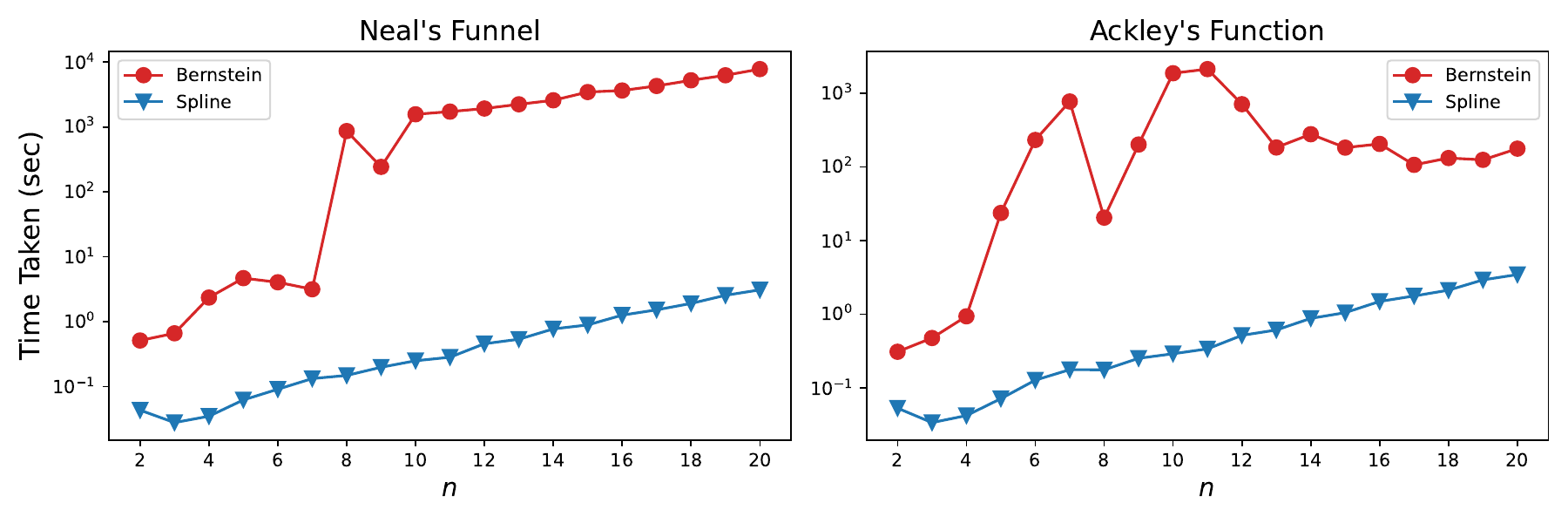}
    \caption{The time taken to fit 2D functions without noise for a Multivariate rational Bernstein Denominator compared with a tensor product smoothing spline with an equivalent degree of freedom.}
    \label{fig:gam_bernstein_time_taken}
\end{figure}

\section{Conclusion}
Rational function approximation has been well-studied for functions without noise and is known to have similar theoretical convergence properties to polynomial approximation. However, few attempts have been made to guarantee pole-free approximations. Poles typically have been treated as a necessary by-product of rational functions, and the preventative measures to avert poles are methods to detect poles \cite{poles_who_cares_2, poles_who_cares} or attempt to stabilize the fit \cite{stablized_sk_algorithm, robust_pade_svd}. 

The AAA algorithm \cite{AAA}, presented by the works of  Nakatsukasa, S\'ete, and Trefethen, represents the current gold-standard method to perform rational function approximation. Indeed, we do not claim that our algorithm can beat the AAA algorithm overall. However, as the original AAA paper writes, `The fact is that the core AAA algorithm risks introducing unwanted poles when applied to problems involving real functions on real intervals'. Moreover, a recent paper published in 2023 by Huybrechs and Trefethen \cite{poles_who_cares} writes, `the appearance
of unwanted poles in AAA approximants is not yet fully understood'. Our experiments show that the unwanted poles can appear even in simple examples.

In this paper, we take an alternative route and put poles at the forefront of our method by guaranteeing rational polynomials with no poles in an interval on the real line. This is performed through Bernstein polynomials and normalized coefficients to force strict positivity in the denominator. This represents a restricted class of rational functions and thus requires extra degrees of freedom to reach the accuracy of traditional rational approximation methods. Yet, the $\mathcal{R}(n,n)$ rational Bernstein denominator numerically displays the expected $\mathcal{P}(2n)$ convergence rate.

Our method typically takes longer than other rational or polynomial approximations, predominantly due to the iterative scheme. However, for our main application in differential equations, this cost pales compared to the typical runtime to numerically solve differential equations using spectral methods, which can take hours or days. Moreover, the compact representation afforded by the rational Bernstein denominator algorithm can dramatically reduce the spectral methods' runtime without any sacrifice to the solution's accuracy.


Paraphrasing the co-creator of the AAA algorithm, Lloyd N. Trefethen, in his textbook on function approximation \cite{approximation_theory_and_practice}, there is no universally best approximation method. The AAA algorithm is best when one requires accurate estimations on data free of noise. We are not trying to replace it. However, when one wants to guarantee that no poles are inside the approximation interval reliably, we believe the Bernstein denominator algorithm provides a promising and robust approach to using rational polynomials on data. 

For further research, we would like to investigate theoretical guarantees on its convergence rate and explore more effective methods for introducing smoothing in rational function approximations. Bernstein polynomials can also be represented in Barycentric forms, leading to more stable ways of computing the rational function and potentially better approximations.

\section*{Acknowledgments}
We thank Johnny Myungwon Lee for the helpful discussions on penalized splines. We also thank the reviewers for their helpful comments and suggested experiments that improved the quality of this paper. 

\bibliographystyle{unsrt}

\end{document}